\documentclass[12pt,a4paper]{article}

\usepackage{theorem,enumerate}
\usepackage{amsmath,latexsym,amssymb,amsfonts, epsfig}
\usepackage{eucal}
\usepackage{mathrsfs}
\usepackage{array}
\usepackage{epstopdf} 

\newcounter{Scounter}
\theorembodyfont{\normalfont\slshape}
\newtheorem{thm}{Theorem}

\newtheorem{cor}[thm]{Corollary}

\newcommand{\Proof}{\medbreak\noindent\textbf{Proof.}\quad}

\newcommand{\qed}{{$\quad\square$\vs{3.6}}}

\newcommand{\vs}[1]{\vspace*{#1 mm}}

\bfseries\normalfont

\makeatletter
\def\thanks#1{%
   \footnotemark
   \edef\@tempa{\noexpand\noexpand\noexpand\footnotetext[\the\c@footnote]}%
   \toks@\expandafter{\@thanks}%
   \toks\tw@{{#1}}
   \xdef\@thanks{\the\toks@\@tempa\the\toks\tw@}}
\makeatother

\begin{document}

\title{A Note on $4$-colorings of Quadrangulations}

\author{
Arthur Hoffmann-Ostenhof\thanks{Technical University of Vienna,
Austria.
Email: {\tt arthurzorroo@gmx.at}},
Atsuhiro Nakamoto \thanks{Yokohama National University, Japan.
Email: {\tt nakamoto@ynu.ac.jp}}
}
\date{}
\maketitle

\begin{abstract}
Let $G$ be a quadrangulation on an orientable surface and let $g$ be a proper vertex-$4$-coloring of $G$. 
A face $F$ of $G$ is said to be a {\em rainbow-face\/} if all four distinct colors appear on its boundary. 
A {\em $(c_1,c_2,c_3,c_4)$-face\/} in $G$ is a rainbow face with colors $c_i$, $i=1,2,3,4$ on the boundary in clockwise order.
We show that the number of $(c_1,c_2,c_3,c_4)$-faces in $G$ equals the number of $(c_4,c_3,c_2,c_1)$-faces. This implies in particular that the number of rainbow-faces of $G$ is even.

\end{abstract}


\section{Introduction}

All graphs considered are loopless and finite.
A \textit{cycle} is a connected $2$-regular graph. 
A \textit{quadrangulation\/} on an oriented surface (i.e a sphere with $k \geq 0$ attached handles) is
a $2$-edge connected graph embedded on the surface 
such that every face boundary is a cycle of length four.
A vertex-coloring is said to be \textit{proper} if any two adjacent vertices receive different colors. 
For standard graph-theoretic terms not defined here, we refer to \cite{Bo, West}.

Let $G$ be a quadrangulation on an orientable surface with a proper vertex-coloring $g: V(G) \mapsto \{1,2,3,4\}$.
Then for each face of $G$, two, three or four different colors appear on its boundary. 
We consider only those faces of $G$ which have all four colors on the boundary and call them \textit{rainbow-faces}.
A rainbow-face $F$ is said to be a \textit{$(c_1,c_2,c_3,c_4)$-face}, where $\{c_1,c_2,c_3,c_4\}=\{1,2,3,4\}$,
if the four distinct colors $c_1,c_2,c_3,c_4$ appear on the boundary of $F$ in the clockwise order,
where the clockwise orientation can be fixed in all faces of $G$ simultaneously, since the surface is orientable.
A $(c_1,c_2,c_3,c_4)$-face is also an $x$-face, where $x \in \{(c_2,c_3,c_4,c_1),(c_3,c_4,c_1,c_2),(c_4,c_1,c_2,c_3)\}$.

It was shown in \cite{Mo} that the number of $(c_1,c_2,c_3,c_4)$-faces
and the number of $(c_4,c_3,c_2,c_1)$-faces have the same parity 
if $G$ is a quadrangulation on the sphere.
Then it was proven in \cite{Hoff} that these two numbers are even equal. 
In this paper, we show that this holds for all orientable surfaces by a very simple proof.
(For example, Figure 1 shows a quadrangulation on the torus, where we identify the top and the bottom,
and the right and left sides, respectively, in the rectangle.)
The idea for the proof comes from the paper \cite{arch}, which points out the surprising phenomenon that
a certain kind of quadrangulations on a nonorientable surface has no proper vertex-$3$-coloring.

\begin{figure}[htpb]\label{quad}
\centering\epsfig{file=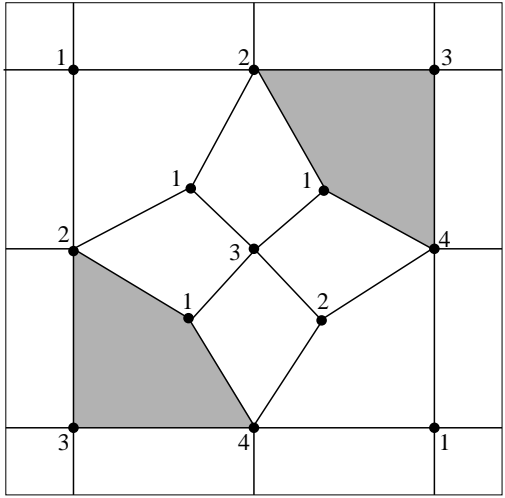,width=2.6 in}
\caption{A quadrangulation of the torus with the $(1,2,3,4)$- and $(4,3,2,1)$-faces in gray.} 
\end{figure}

\section{Main result and proof}

\begin{thm}\label{main}
Let $G$ be a quadrangulation on an orientable surface with a proper vertex-$4$-coloring $g: V(G) \mapsto \{1,2,3,4\}$.
Then the number of $(c_1,c_2,c_3,c_4)$-faces in $G$ equals the number of $(c_4,c_3,c_2,c_1)$-faces,
where $\{c_1,c_2,c_3,c_4\}=\{1,2,3,4\}$.
\end{thm}

\Proof
We may assume that $(c_1,c_2,c_3,c_4)=(1,2,3,4)$ otherwise we permute the four colors in a suitable manner.
We orient all edges of $G$ by giving a direction to each edge $xy$ from $x$ to $y$ if and only if $g(x)>g(y)$. 
Let $F$ be a face of $G$ and let $F^+$ (resp., $F^-$) denote the set of all boundary edges of $F$
which are oriented clockwise (resp., anti-clockwise) in $F$. 
Set $\Delta(F):= |F^+| - |F^-|$.
Let $\mathcal F(G)$ denote the set of all faces of $G$ and let $F ^*(G)$ denote the set of all $(1,2,3,4)$- 
and $(4,3,2,1)$-faces in $G$.
Then,
$$ \sum_{F\in \mathcal F(G)}\Delta (F)= \sum_{F \in \mathcal F(G)-F^*(G)}\Delta(F) + \sum_{F \in F^*(G)}\Delta(F) = 0\,\,$$
since each $e \in E(G)$ is contained in precisely two faces, say $F$ and $F'$, sharing $e$,
in which if $e$ has the clockwise orientation in $F$, then $e$ has the anti-clockwise orientation in $F'$, and vice versa. 
On the other hand, we see that $\Delta(F)=0$ if and only if $F$ is neither a $(1,2,3,4)$-face nor a $(4,3,2,1)$-face in $G$.
Therefore,
since $\sum_{F \in \mathcal F(G)-F^*(G)}\Delta(F)=0$, we have
$$\sum_{F \in F^*(G)}\Delta(F) = 0\,\,.$$ 
Observe that $\Delta(F)=-\Delta(F') \ne 0$ if $F$ is a $(1,2,3,4)$-face and $F'$ is a $(4,3,2,1)$-face,
and hence the number of $(1,2,3,4)$-faces equals the number of $(4,3,2,1)$-faces. \qed

\bigskip


\begin{cor}\label{}
Let $G$ be a quadrangulation of an orientable surface.
Then $G$ admits no proper vertex $4$-coloring of $G$ such that the number of rainbow faces is odd.
\end{cor}

\section*{Acknowledgments}
This work was supported by the Austrian Science Fund (FWF Projects: P 20543, P 26686).


\end{document}